\documentclass[11pt,leqno]{article}

\usepackage{amsfonts,enumerate,amsmath,amsfonts,xspace,array,delarray,theorem,amssymb,epsfig}

\theoremstyle{break}\theorembodyfont{\it}\theoremheaderfont{\normalfont\bfseries}

\newtheorem{theo}{Theorem}

\newtheorem{defi}[theo]{Definition}

\newtheorem{lem}[theo]{Lemma}

\newtheorem{prop}[theo]{Proposition}

\newtheorem{coro}[theo]{Corollary}

\newenvironment{proof}{\noindent{\bf Proof: }}
                {\leavevmode\unskip\nobreak\hskip2em plus1fill
                $\scriptstyle\bullet$\vskip\theorempostskipamount\par}

\def\hor{H\"ormander\xspace}
\def\LP{Littlewood--Paley\xspace}

\def\freccia{\longrightarrow}

\def\R{{\mathbb R}}

\def\N{{\mathbb N}}
\def\Z{{\mathbb Z}}
\def\LL{{\cal L}}
\def\SS{{\cal S}}
\def\rn{\R^n}
\def\H{{\mathbb H}}

\def\dg{{\cal D}(G)}
\def\dig{{\cal D'}(G)}
\def\sg{{\cal S}(G)}
\def\sig{{\cal S}'(G)}
\def\lp#1{L^{#1}(G)}

\def\pair#1#2{\langle#1,#2\rangle}

\def\mez{\frac{1}{2}}
\def\lm{\sqrt{\lambda}}
\def\htm{h_t}
\def\wht{\widehat{\htm}}
\def\lmez{\sqrt{\LL}}

\def\cb#1#2{\left(\begin{array}{c}#1\\#2\end{array}\right)}

\def\beq{\begin{eqnarray*}}
\def\eeq{\end{eqnarray*}}
\let\lt=<
\let\gt=>

\begin{document}
\title{\LP decompositions and Besov spaces on Lie groups of polynomial growth}
\date{}
\author{\large {Giulia Furioli, Camillo Melzi and Alessandro Veneruso}}

\maketitle
\footnotetext{The authors were partially supported by GNAMPA - Progetto
``{\em Calcolo funzionale
per generatori di semigruppi ed analisi armonica su gruppi}'', 2003.}

\begin{abstract}
We introduce a Littlewood--Paley decomposition related to any
sub-Laplacian on
a Lie group $G$ of polynomial volume growth;  this allows us to prove
a  Littlewood--Paley theorem in this general setting and
to provide a dyadic
characterization of Besov spaces $B^{s,q}_p(G)$, $s\in\R$, equivalent to
the classical definition through the heat kernel.
\end{abstract}

\section{Introduction}\label{intro}\leavevmode\par
\LP decompositions are a powerful tool in investigating
deep properties of function spaces of distributions.
Let us recall briefly the classical  construction in $\R^n$.
Let  $\varphi \in C^\infty
(\R)$ be an even function  such that $0\leq \varphi \leq 1$, $\varphi = 1$
in $[0,\frac 14]$ and $\varphi =0$ in
$[1,\infty).$ Let $\psi(\lambda)= \varphi (\frac \lambda 4)
-\varphi(\lambda)$, so that
$\rm {supp\ } \psi \subset \{  \frac 14 \leq |\lambda| \leq 4\}$.
We have the following partition of unity on the frequency space of the
Fourier transform:
\[
1= \varphi (|\xi|^2) + \sum_{j= 0}^\infty \psi(2^{-2j}|\xi|^2), \quad
\xi\in\R^n.
\]
This gives the identity in $\SS'(\R^n)$:
\[
\hat u= \varphi (|\cdot|^2) \hat u+ \sum_{j= 0}^\infty
\psi(2^{-2j}|\cdot|^2)\hat u, \quad u\in \SS'(\R^n)
\]
and, denoting by $S_0u$ e $\Delta_j u$ respectively
\[
\widehat {S_0u} =\varphi (|\cdot|^2)\hat u, \quad \widehat{\Delta_j u}
=\psi(2^{-2j}|\cdot|^2)\hat u,
\]
we obtain the \LP decomposition in $\SS'(\R^n):$
\begin{equation}\label{decomp}
u= S_0u+ \sum_{j= 0}^\infty  \Delta_ju,\quad u\in \SS'(\R^n).
\end{equation}
The following fundamental theorem holds (see e.g. \cite{Ste}):
\begin{theo}[\LP]
Let $1<p<\infty$ and $u\in\SS'(\R^n)$. Then $u\in L^p(\rn) $ if and only
if $S_0 u\in L^p(\rn) $ and $\left( \sum_{j= 0}^\infty
|\Delta_ju|^2\right)^{\frac 12} \in
L^p(\rn) $. Moreover there exists a constant $C_p>1$, which depends only on
$p$, such that
\[
C_p^{-1}\|u\|_{L^p (\rn)}\leq\|S_0 u\|_{L^p(\rn)}+\left \|
\left( \sum_{j= 0}^\infty  |\Delta_ju|^2\right )^{\frac 12}
\right \|_{L^p(\rn)}\leq
C_p\|u\|_{L^p (\rn)},\qquad u\in L^p(\rn).
\]
\end{theo}
The proof of this theorem is based on the classical H\"ormander--Mihlin
$L^p$-multiplier theorem (\cite{Hor}) and on the uniform estimates for the
norms of the convolution operators $\Delta_j $ on $L^p(\rn)$.
The main purpose of this paper is to prove a \LP theorem on Lie groups of
polynomial volume growth, with respect to any sub-Laplacian.

In view of extending the previous construction to a general Lie group of
polynomial growth it is more convenient to see the decomposition
\eqref{decomp} in terms of multipliers of the Laplacian $\Delta =
-\sum_{i=1}^n \frac{\partial^2}{\partial x_i^2}$. Starting from the
spectral decomposition of the Laplacian
\[
\Delta = \int_0^\infty \lambda dE_\lambda
\]
and from the functions $\varphi$, $\psi\in L^\infty(\rn)$ previously
introduced, we consider the multiplier operators
\beq
&&\varphi (\Delta)  = \int_0^\infty \varphi(\lambda)dE_\lambda \, \\
&&\psi(2^{-2j}\Delta) = \int_0^\infty \psi(2^{-2j}\lambda)dE_\lambda.
\eeq
So we have the identifications between operators
\[
\varphi (\Delta) f =S_0f, \quad \psi(2^{-2j}\Delta) f= \Delta_j f, \quad
f\in L^2 (\rn).
\]
If we denote by $\Psi_j$ the convolution kernel of the operator
$\psi(2^{-2j}\Delta)$, due to the dilation structure of $\rn$ we have the
scaling formula
\begin{equation}\label{dilataz}
\Psi_j(x)= 2^{nj}\Psi_0(2^j x),\quad x\in \R^n, j\in \N
\end{equation}
that allows us to easily obtain the estimates for the norms of the
operators $\Delta_j$ only from the operator $\Delta_0$.
If we now consider a stratified Lie group $G$, endowed with its natural
dilation structure, and if
$\Delta$ is a sub-Laplacian on $G$ invariant with respect to the family of
dilations, the scaling
formula
\eqref{dilataz} still holds, where $n$ is the homogeneous dimension of $G$.
In this case the way to
prove a \LP theorem through a \hor--Mihlin multiplier theorem is based on
classical
techniques (\cite{FoSt},
\cite{DeMa}, \cite{MaMe}). In the particular case of the Heisenberg group
$\H_n$, we can deduce a
\LP theorem also for the full Laplacian from the results by M\"uller, Ricci
and Stein \cite{MRS1},
\cite{MRS2}; but their techniques, based on the Fourier transform on
$\H_n$, do not fit to the case
of Lie groups of polynomial growth. Alexopoulos in \cite{Ale2} proved a
\hor--Mihlin multiplier
theorem for any sub-Laplacian in the general setting of Lie groups of
polynomial growth;
nevertheless, such result does not provide directly the uniform estimates
for the norms of the
operators $\Delta_j$ we need to deduce a \LP theorem.

The main result of this paper is Proposition \ref{main}, which allows us to
deduce the uniform estimates for the norms of the operators $\Delta_j$ and
to prove a \LP decomposition in $\SS'(G)$ (Proposition \ref{lpd}) and a \LP
theorem (Theorem \ref{llpp}).
As an application of such decomposition, we finally provide a dyadic
characterization of Besov
spaces $B^{s,q}_{p}(G),$ $s\in \mathbb R$, equivalent to the classical
definition through the heat
kernel (Proposition \ref{besov-eq}).\par\medskip
We would like to thank Stefano Meda for several helpful discussions.

\section{Notation and preliminaries}\label{Notpre}\leavevmode\par
In this paper
$\N$ denotes the set of nonnegative integers, $\Z_{+}$ the set of positive
integers and
$\R_{+}$ the set of positive real numbers. For $p\in[1,\infty]$ we denote
by $p'$ the
conjugate index of $p$, such that
$\frac{1}{p}+\frac{1}{p'}=1$.\par\smallskip In this section we
recall some basic facts about Lie groups of polynomial growth. For the
proofs and further
information, see e.g.\
\cite{VSC} and the references given therein.\par\smallskip  Let $G$ be a
connected Lie group, and
let us fix a left-invariant Haar measure
$dx$ on $G$. We will denote by $|A|$ the measure of a
measurable subset $A$ of $G$ and by $\chi_A$ its characteristic
function.\par\smallskip
We assume that
$G$ has polynomial volume growth, i.e., if
$U$ is a compact neighbourhood of the identity element $e$ of $G$, then
there is a constant $C>0$
such that
$|U^n|\leq Cn^C$, $n\in\Z_{+}$. Then $G$ is unimodular. Furthermore, there
exists $D\in\N$, which
does not depend on
$U$, such that
\begin{equation}|U^n|\sim n^D\qquad {\rm for}\
n\rightarrow\infty\label{un}.\end{equation} For instance, every connected
nilpotent Lie group has
polynomial volume growth.\par\smallskip
The convolution of two functions $f$ and $g$ on $G$ is
defined by
\[f*g(x)=\int_G f(y)g(y^{-1}x)\,dy,\qquad x\in G\]
and satisfies the Young's inequality (where
$1+\frac{1}{r}=\frac{1}{p}+\frac{1}{q}$)
\[\|f*g\|_{\lp r}\leq\|f\|_{\lp p}\|g\|_{\lp q}.\]
The space $\dg$ of test
functions and the space
$\dig$ of distributions are
defined in the usual way (see
\cite{Ehr}). The convolution of $\varphi\in\dg$ and $u\in\dig$ is
defined as usual:
\begin{equation}\pair{\varphi*u}{\psi}=\pair{u}{\check{\varphi}*\psi},\qquad\psi
\in\dg
\label{convoluz}\end{equation}
where
\[\check{\varphi}(x)=\varphi(x^{-1}),\qquad x\in G.\]\par\smallskip
Let $X_1,\ldots,X_k$ be left-invariant vector fields on
$G$ which satisfy the H\"ormander's condition, i.e.\ they generate,
together with their successive
Lie brackets
$[X_{i_1},[\ldots,X_{i_{\alpha}}]\cdots]$, the Lie algebra of $G$. For
$I=(i_1,\ldots,i_{\beta})\in\{1,\ldots,k\}^{\beta}$ ($\beta\in\N$) we put
$|I|=\beta$
and $X^I=X_{i_1}\cdots X_{i_{\beta}}$, with the convention that $X^I=id$ if
$\beta=0$.\par\smallskip
To $X_1,\ldots,X_k$
is associated, in a canonical way, the control distance
$\rho$, which is left-invariant and compatible with the topology on $G$.
For any $x\in G$ we put
$|x|=\rho(e,x)$. The properties of $\rho$ imply that $|xy|\leq|x|+|y|$ for
any $x,y\in
G$. Furthermore, for any
$r>0$ we put
$V(r)=|B(e,r)|$ where $B(e,r)=\{x\in G:|x|<r\}$. By~(\ref{un}) we have
\[V(r)\sim r^D\qquad {\rm for}\
r\rightarrow\infty.\]
On the other hand, there exists $d\in\N$ such that
\[V(r)\sim r^d\qquad {\rm for}\
r\rightarrow 0.\]
These estimates imply the ``doubling property'': there exists $K>0$ such that
\begin{equation} V(2r)\leq KV(r),\qquad
r>0.\label{doubling}\end{equation}\par\smallskip
We consider the
sub-Laplacian
\[\LL=-\sum_{j=1}^k X_j^2\]
which is a positive self-adjoint operator, having as domain of
definition the space of all functions $f\in\lp 2$ such that $\LL f\in\lp 2$.
So, by the spectral theorem, for any bounded Borel function $m$ on
$[0,\infty)$ we can define the
operator $m(\LL)$ which is bounded on $\lp 2$. Since the point 0 may be
neglected in
the spectral resolution of $\LL$ (see \cite{Chr}, \cite{Ale2}), we consider
that
the function
$m$ is defined on $\R_{+}$. Furthermore, the operator
$m(\LL)$ admits a kernel
$M\in\dig$ which satisfies $m(\LL)f=f*M$ for any $f\in\dg$. We recall the
following
well-known results:
\begin{theo}[\cite{Ale2}] Put $N=1+\max\{[\frac{d}{2}],[\frac{D}{2}]\}$. If
$m\in C^N(\R_{+})$ and
$\sup_{\lambda>0} \lambda^r|m^{(r)}(\lambda)|<\infty$ for any
$r\in\{0,\ldots,N\}$, then $m(\LL)$
extends to a bounded operator on $\lp p$,
$1<p<\infty$.\label{alexop}\end{theo}
\begin{prop} Let $\{m_n\}_{n\in\N}$ be a sequence of bounded Borel
functions on $\R_{+}$ which
converges at every point to a bounded Borel function $m$. Suppose also that the
sequence $\{\|m_n\|_{L^{\infty}(\R_{+})}\}_{n\in\N}$ is bounded. Then the
sequence
$\{M_n\}_{n\in\N}$, where $M_n$ is the kernel of the operator $m_n(\LL)$,
converges in $\dig$ to
the kernel $M$ of the operator $m(\LL)$.
\label{dixmal}\end{prop}
\begin{proof} By the spectral theorem $m_n(\LL)f\rightarrow m(\LL)f$ in
$\lp 2$ for $n\rightarrow
\infty$ for every
$f\in\lp 2$. In particular, $f*M_n\rightarrow f*M$ in $\dig$ for $n\rightarrow
\infty$ for every
$f\in\dg$. Fix $\varphi\in\dg$. By \cite[Th\'eor\`eme 3.1]{DiMa} the
function $\varphi$ can be
written as a finite sum $\varphi=\sum_{j=1}^r\psi_j*\chi_j$ with
$\psi_j,\chi_j$ in $\dg$. So by
(\ref{convoluz})
\[\pair{M_n}{\varphi}=\sum_{j=1}^r\pair{\check{\psi_j}*M_n}
{\chi_j}\stackrel{n\rightarrow
\infty}{\longrightarrow}\sum_{j=1}^r\pair{\check{\psi_j}*M}{\chi_j}=
\pair{M}{\varphi}.\]
\end{proof}\par\smallskip
We introduce the Schwartz space $\sg$ and its dual space $\sig$ as
in
\cite{Skr1},\cite{Skr2}. The definition does not depend on
$X_1,\ldots,X_k$. However, a useful characterization of $\sg$ is the
following: a function
$f\in C^{\infty}(G)$ is in
$\sg$ if and only if all the seminorms
\[p_{\alpha,I}(f)=\sup_{x\in G}(1+|x|)^{\alpha}|X^I f(x)|\]
($\alpha\in\N,\ I\in\bigcup_{\beta\in\N}\{1,\ldots,k\}^{\beta}$) are
finite. The space $\sg$
endowed with this family of seminorms is a Fr\'echet space. It is easy to
show that
$\sg*\sg\subset\sg\subset\lp p$ for $1\leq p\leq\infty$.\par\smallskip
The heat kernel $p_t$, i.e.\ the kernel of the operator
$e^{-t\LL}$ ($t>0$), is a positive $C^{\infty}$ function which satisfies
$\int_G p_t(x)\,dx=1$.
Moreover, for any
$I\in\bigcup_{\beta\in\N}\{1,\ldots,k\}^{\beta}$ there exists
$C>0$ such that the following estimates hold:
\begin{eqnarray}p_t(x)&\leq&CV(\sqrt t)^{-1}e^{-\frac{|x|^2}{Ct}},\qquad
x\in G,\
t>0;\label{ptx}\\
|X^I p_t(x)|&\leq&Ct^{-\frac{d+|I|}{2}}e^{-\frac{|x|^2}{Ct}},\qquad x\in G,\
0<t\leq 1.\label{xiptx}\end{eqnarray}
In particular, estimate (\ref{xiptx}) implies that $p_t\in\sg$ for
$t\in(0,1]$. Since
$p_{t_1+t_2}=p_{t_1}*p_{t_2}$ for any $t_1,t_2>0$, it follows that
$p_t\in\sg$ for any
$t>0$. Furthermore, estimates (\ref{ptx}) and (\ref{xiptx}) yield the following
\begin{prop} For any $\alpha\in\N$,
$I\in\bigcup_{\beta\in\N}\{1,\ldots,k\}^{\beta}$ and
$p\in[1,\infty]$ there exists
$C>0$ such that the following estimates hold:
\begin{eqnarray}\|(1+|\cdot|)^{\alpha}p_t(\cdot)\|_{\lp
p}&\leq&C(1+\sqrt t)^{\alpha}V(\sqrt t)^{-\frac{1}{p'}},\qquad
t>0;\label{ptp}\\
\|(1+|\cdot|)^{\alpha}X^I p_t(\cdot)\|_{\lp
p}&\leq&Ct^{-(\frac{d}{2p'}+\frac{|I|}{2})},\qquad
0<t\leq 1.\label{xiptp}\end{eqnarray}\label{nuova}\end{prop}
\begin{proof} In this proof we will denote by
$C$ a positive constant which will not be
necessarily the same at each occurrence,
with the convention that $C$ can depend only on $G$ and on
$\alpha,I,p$.\par\smallskip
Fix $t>0$. First we note that
\begin{equation}\sup_{\rho\geq
0}(1+\rho)^{\alpha}e^{-\frac{\rho^2}{Ct}}\leq
C(1+\sqrt t)^{\alpha}\label{suprho}\end{equation} as is easy to verify by
calculating the
maximum of the function $\rho\mapsto
(1+\rho)^{\alpha}e^{-\frac{\rho^2}{Ct}}$ in $[0,\infty)$. For
$p=\infty$ estimates (\ref{ptp}) and (\ref{xiptp}) follow directly by
(\ref{suprho}) and by
(\ref{ptx}) and (\ref{xiptx}), respectively. For $1\leq p<\infty$ we use
the fact that
\begin{equation}\int_G e^{-\frac{|x|^2}{Ct}}\,dx\leq CV(\sqrt
t)\label{inte}\end{equation}
(see \cite[page 111]{VSC}). So by (\ref{ptx}), (\ref{suprho}) and (\ref{inte})
\beq\|(1+|\cdot|)^{\alpha}p_t(\cdot)\|_{\lp p}&\leq&CV(\sqrt
t)^{-1}\left(\int_G(1+|x|)^{\alpha
p}e^{-\frac{p|x|^2}{Ct}}\,dx\right)^{\frac{1}{p}}\\
&\leq&CV(\sqrt t)^{-1}\left(\sup_{x\in G}(1+|x|)^{\alpha
p}e^{-\frac{p|x|^2}{2Ct}}\right)^{\frac{1}{p}}\left(\int_G
e^{-\frac{p|x|^2}{2Ct}}\,dx\right)^{\frac{1}{p}}\\ &\leq&C(1+\sqrt
t)^{\alpha}V(\sqrt
t)^{-\frac{1}{p'}}.\eeq
So (\ref{ptp}) is proved. The proof of (\ref{xiptp}) is
analogous: it only uses (\ref{xiptx}) instead of (\ref{ptx}).
\end{proof}

\section{Weighted estimates and Schwartz kernels}\label{Weiest}\leavevmode\par
In this section we prove some weighted
$L^p$ estimates for the kernel $M_t$ of the operator $m(t\LL)$, $t>0$, when
$m$ is a sufficiently regular function on $\R_{+}$. In the  case where
$G$ is a stratified group and
$\LL$ is invariant with respect to the dilation structure of $G$, the
relation between $M_t$ and
$M_1$ stated in \cite[Lemma 6.29]{FoSt} makes such estimates much easier to
prove, since they
follow directly from the corresponding estimates for $M_1$
(\cite[Lemmas 6.35 and 6.36]{FoSt}, \cite[Lemma 3.1]{DeMa}, \cite[Lemma
1.2]{MaMe}). In the general
case the situation is more complicated and we have to use more
sophisticated techniques due
essentially to Alexopoulos. In particular, we will use the following lemma (see
\cite{Ale1}, \cite{Ale2} and the references given there):
\begin{lem} Fix $\delta>0$, $n\in\Z_{+}$ and $h\in C^n(\R)$ with compact
support.
Then there is an even function $g\in L^1(\R)\cap L^{\infty}(\R)\cap C(\R)$
such that
{\rm supp}$\,\hat g\subset[-\delta,\delta]$ and
\[\sup_{\lambda\in\R}|h(\lambda)-h*g(\lambda)|\leq
C\delta^{-n}\sup_{\sigma\in\R}|h^{(n)}(\sigma)|\]
where $C$ is a positive constant which depends only on
$n$ but not on $\delta,h,g$.\label{alexo}\end{lem}\par\smallskip
Throughout this section we will use the following notation: for any
$n\in\N$ and for any
$m\in C^n(\R_{+})$ we put
\[\|m\|_{(n)}=\sup_{\stackrel{\scriptstyle 0\leq
r\leq n}{\lambda>0}}(1+\lambda)^n|m^{(r)}(\lambda)|.\]
Moreover, the constants $d$ and $K$ are the same which have been introduced
in Section~\ref{Notpre}.
\begin{prop} Fix $\alpha\in\N$,
$I\in\bigcup_{\beta\in\N}\{1,\ldots,k\}^{\beta}$ and
$p\in[1,\infty]$. There exist $C>0$ and
$n\in\Z_{+}$ such that for any $m\in C^n(\R_{+})$ with $\|m\|_{(n)}<\infty$
the kernel $M_t$ of the operator
$m(t\LL)$, $t>0$, satisfies the following estimates:
\begin{eqnarray}\|(1+|\cdot|)^{\alpha}M_t(\cdot)\|_{\lp p}&\leq&
C(1+\sqrt t)^{\alpha}V(\sqrt t)^{-\frac{1}{p'}}\|m\|_{(n)},\qquad
t>0;\label{normelp}\\
\|(1+|\cdot|)^{\alpha}X^I M_t(\cdot)\|_{\lp p}&\leq&
Ct^{-(\frac{d}{2p'}+\frac{|I|}{2})}\|m\|_{(n)},\qquad 0<t\leq
1.\label{aggiunta}\end{eqnarray}
\label{main}\end{prop}\par\smallskip
{\em Remark:}\quad The case where $\alpha=|I|=0$ and $p=1$ is particularly
interesting: it simply
reads
\[\|M_t\|_{\lp 1}\leq C\|m\|_{(n)},\qquad t>0.\]\par\smallskip
\begin{proof}  In this proof we will denote by
$C$ a positive constant which will not be
necessarily the same at each occurrence,
with the convention that $C$ can depend only on $G$ and on
$\alpha,I,p$.
The proof consists of some
steps.\par\smallskip
{\bf Step 1.}\quad We prove (\ref{normelp}) for
$p=1$, with the additional assumption that $m=0$ in
$[2,\infty)$. Fix $t>0$ and fix $n\in\Z_{+}$ which will be choosen later.
We consider the function $h_t$ on $\R$ defined by
\[h_t(\sigma)=e^{t\sigma^2}m(t\sigma^2),\qquad\sigma\in\R.\]
By the assumptions on $m$ we have
\begin{equation}\|h_t\|_{L^{\infty}(\R)}\leq
e^2\|m\|_{(0)}.\label{qualcosa}\end{equation}
Moreover
\[m(t\lambda)=e^{-t\lambda}h_t(\sqrt{\lambda}),\qquad\lambda>0.\]
So by the spectral theorem
\begin{equation} M_t=h_t(\sqrt{\LL})p_t\in\lp 2\label{mt}\end{equation}
and
\beq\|M_t\|_{\lp 2}&\leq&\|h_t\|_{L^{\infty}(\R)}\|p_t\|_{\lp 2}\\
&\leq&C V(\sqrt t)^{-\frac{1}{2}}\|m\|_{(0)}\eeq
by (\ref{ptp}) and (\ref{qualcosa}). Then
\begin{eqnarray}\int_{|x|<\sqrt t}(1+|x|)^{\alpha}|M_t(x)|\,dx&\leq&
\left(\int_{|x|<\sqrt t}(1+|x|)^{2\alpha}\,dx\right)^{\mez}
\left(\int_{|x|<\sqrt t}|M_t(x)|^2\,dx\right)^{\mez}\nonumber\\
&\leq&V(\sqrt t)^{\mez}(1+\sqrt t)^{\alpha}\|M_t\|_{\lp 2}\nonumber\\
&\leq&C (1+\sqrt t)^{\alpha}\|m\|_{(0)}.
\label{civuole}\end{eqnarray}
On the other hand it follows from (\ref{mt}) that
\begin{equation}\int_{|x|\geq
\sqrt t}(1+|x|)^{\alpha}|M_t(x)|\,dx=\sum_{j=0}^{\infty}
\left(\int_{A_{t,j}}(1+|x|)^{\alpha}|M^{(1)}_{t,j}(x)|\,dx+
\int_{A_{t,j}}(1+|x|)^{\alpha}|M^{(2)}_{t,j}(x)|\,dx\right)\label{spezza}\end{equation}
where:
\beq A_{t,j}&=&\{x\in G: 2^j\sqrt t\leq|x|<2^{j+1}\sqrt t\};\\
M^{(1)}_{t,j}&=&h_t(\sqrt{\LL})(p_t\chi_{\{y\in
G:|y|<2^{j-1}\sqrt t\}});\\
M^{(2)}_{t,j}&=&h_t(\sqrt{\LL})(p_t\chi_{\{y\in
G:|y|\geq 2^{j-1}\sqrt t\}}).\eeq
For every $j\in\N$ and for $i=1,2$ we have
\beq\int_{A_{t,j}}(1+|x|)^{\alpha}|M^{(i)}_{t,j}(x)|\,dx&\leq&
\left(\int_{A_{t,j}}(1+|x|)^{2\alpha}\,dx\right)^{\mez}
\left(\int_{A_{t,j}}|M^{(i)}_{t,j}(x)|^2\,dx\right)^{\mez}\\
&\leq&V(2^{j+1}\sqrt t)^{\mez}(1+2^{j+1}\sqrt
t)^{\alpha}\|M^{(i)}_{t,j}\|_{L^2(A_{t,j})}.\eeq
The fact that
$(1+2^{j+1}\sqrt t)^{\alpha}\leq C 2^{j\alpha}(1+\sqrt t)^{\alpha}$
and the doubling
property (\ref{doubling}) imply
\begin{equation}\int_{A_{t,j}}(1+|x|)^{\alpha}|M^{(i)}_{t,j}(x)|\,dx\leq
CK^{\frac{j}{2}}
V(\sqrt
t)^{\mez}2^{j\alpha}(1+\sqrt t)^{\alpha}\|M^{(i)}_{t,j}\|_{L^2(A_{t,j})}.
\label{iaj}\end{equation}
In order to estimate $\|M^{(1)}_{t,j}\|_{L^2(A_{t,j})}$, first of all we
suppose $n\geq 2$ and
we note that
\[h_t(\sqrt{\lambda})=\frac{1}{2\pi}\int_{\R}\wht(s)\cos(s\lm)\,ds,\qquad\lambda
>0\]
since $\htm$ is an even function in $L^1(\R)$ whose
Fourier transform is in $L^1(\R)$. Then for a.e.\ $x\in A_{t,j}$
\begin{equation}M^{(1)}_{t,j}(x)=\frac{1}{2\pi}\int_{\R}\wht(s)\left(\cos(s\lmez
)(p_t\chi_{\{y\in
G:|y|<2^{j-1}\sqrt t\}})\right)(x)\,ds.\label{htmlpt}\end{equation}
Now we use the fact that the kernel $G_s$ of the operator $\cos(s\lmez)$,
$s\in\R$, satisfies the
following property (see \cite{Mel}):
\[{\rm supp}\,G_s\subset\{y\in G:|y|\leq|s|\}.\]
So, for $x\in A_{t,j}$ and $|s|\leq 2^{j-1}\sqrt t$ we have
\begin{equation}\left(\cos(s\lmez)(p_t\chi_{\{y\in
G:|y|<2^{j-1}\sqrt t\}})\right)(x)=0.\label{fps}\end{equation} By
Lemma~\ref{alexo} we can take an even function
$g_{t,j}\in L^1(\R)\cap L^{\infty}(\R)\cap C(\R)$ such that
{\rm supp}$\,\widehat{g_{t,j}}\subset[-2^{j-1}\sqrt t,2^{j-1}\sqrt t]$ and
\begin{equation}\sup_{\lambda\in\R}|\htm(\lambda)-\htm*g_{t,j}(\lambda)|\leq
C
2^{-jn}t^{-\frac{n}{2}}\sup_{\sigma\in\R}|\htm^{(n)}(\sigma)|.\label{suptr}\end{equation}
The support property of $\widehat{g_{t,j}}$ and property (\ref{fps}) imply
that for a.e.\ $x\in
A_{t,j}$
\[\int_{\R}\wht(s)\widehat{g_{t,j}}(s)\left(\cos(s\lmez)(p_t\chi_{\{y\in
G:|y|<2^{j-1}\sqrt t\}})\right)(x)\,ds=0.\]
So formula (\ref{htmlpt}) can be rewritten in the following way: for a.e.\
$x\in A_{t,j}$
\beq M^{(1)}_{t,j}(x)&=&
\frac{1}{2\pi}\int_{\R}\left(\wht(s)-\wht(s)\widehat{g_{t,j}}(s)\right)
\left(\cos(s\lmez)(p_t\chi_{\{y\in G:|y|<2^{j-1}\sqrt t\}})\right)(x)\,ds\\
&=&\left((h_t-\htm*g_{t,j})(\lmez)(p_t\chi_{\{y\in
G:|y|<2^{j-1}\sqrt t\}})\right)(x)\eeq
since $\htm-\htm*g_{t,j}$ is an even function in $L^1(\R)\cap
L^{\infty}(\R)\cap C(\R)$ whose
Fourier transform is in $L^1(\R)$. Then
\beq\|M^{(1)}_{t,j}\|_{L^2(A_{t,j})}&\leq&
\left\|(h_t-\htm*g_{t,j})(\lmez)(p_t\chi_{\{y\in
G:|y|<2^{j-1}\sqrt t\}})\right\|_{L^2(G)}\\
&\leq&\|h_t-\htm*g_{t,j}\|_{L^{\infty}(\R)}\|p_t\chi_{\{y\in
G:|y|<2^{j-1}\sqrt t\}}\|_{\lp 2}.\eeq
We apply estimate (\ref{suptr}) to the first factor and we estimate the
second factor simply by
(\ref{ptx}) and (\ref{doubling}), so that
\beq\|p_t\chi_{\{y\in
G:|y|<2^{j-1}\sqrt t\}}\|_{\lp 2}&\leq&C V(\sqrt t)^{-1}V(2^{j-1}\sqrt
t)^{\mez}\\
&\leq&CK^{\frac{j}{2}} V(\sqrt t)^{-\mez}\eeq
and then
\begin{equation}\|M^{(1)}_{t,j}\|_{L^2(A_{t,j})}\leq C
2^{-jn}t^{-\frac{n}{2}}K^{\frac{j}{2}}
V(\sqrt
t)^{-\mez}\|\htm^{(n)}\|_{L^{\infty}(\R)}.\label{cassetto}\end{equation}
We still have to estimate $\|\htm^{(n)}\|_{L^{\infty}(\R)}$. Note that
\[\htm^{(n)}(\sigma)=t^{\frac{n}{2}}
h_1^{(n)}(t^{\mez}\sigma),\qquad\sigma\in\R\]
and
\[h_1^{(n)}(\lambda)=
\sum_{r=0}^n\cb{n}{r}\left(\frac{d^{n-r}}{d\lambda^{n-r}}\right)(e^{\lambda^2})
\left(\frac{d^r}{d\lambda^r}\right)(m(\lambda^2)),\qquad\lambda\in\R.\]
So
\begin{equation}\|\htm^{(n)}\|_{L^{\infty}(\R)}=t^{\frac{n}{2}}\|h_1^{(n)}\|_{L^
{\infty}(\R)}\leq
C t^{\frac{n}{2}}\|m\|_{(n)}.\label{htmn}\end{equation}
Inequalities (\ref{cassetto}) and (\ref{htmn}) give
\begin{equation}\|M^{(1)}_{t,j}\|_{L^2(A_{t,j})}\leq C
2^{-jn}K^{\frac{j}{2}}
V(\sqrt t)^{-\mez}\|m\|_{(n)}.\label{cisiamo}\end{equation}
Now we choose $n>\alpha+\log_2 K$. Then estimates (\ref{iaj}) and
(\ref{cisiamo}) yield
\begin{equation}\sum_{j=0}^{\infty}
\left(\int_{A_{t,j}}(1+|x|)^{\alpha}|M^{(1)}_{t,j}(x)|\,dx\right)\leq
C(1+\sqrt t)^{\alpha}\|m\|_{(n)}.\label{muno}\end{equation}
In order to estimate $\|M^{(2)}_{t,j}\|_{L^2(A_{t,j})}$, first of all we
use the properties of $p_t$
to prove that
\beq\int_{|y|\geq 2^{j-1}\sqrt t}(p_t(y))^2\,dy&\leq&\left(\sup_{|w|\geq
2^{j-1}\sqrt t}p_t(w)\right)
\int_G p_t(y)\,dy\\
&\leq&C V(\sqrt t)^{-1}e^{-\frac{2^{2j}}{C}}\eeq
and then
\beq\|M^{(2)}_{t,j}\|_{L^2(A_{t,j})}&\leq&\|M^{(2)}_{t,j}\|_{\lp 2}\\
&\leq&C\|p_t\chi_{\{y\in
G:|y|\geq 2^{j-1}\sqrt t\}}\|_{\lp 2}\|m\|_{(0)}\\
&\leq&CV(\sqrt t)^{-\mez}e^{-\frac{2^{2j}}{C}}\|m\|_{(0)}.\eeq
This estimate and estimate (\ref{iaj}) yield
\begin{equation}\sum_{j=0}^{\infty}
\left(\int_{A_{t,j}}(1+|x|)^{\alpha}|M^{(2)}_{t,j}(x)|\,dx\right)\leq
C(1+\sqrt t)^{\alpha}\|m\|_{(0)}.\label{mdue}\end{equation}
So both terms of the right-hand side of (\ref{spezza}) can be estimated by
(\ref{muno}) and
(\ref{mdue}), respectively. Thus, taking into account also (\ref{civuole}),
we obtain
\begin{equation}\int_G(1+|x|)^{\alpha}|M_t(x)|\,dx\leq
C(1+\sqrt t)^{\alpha}\|m\|_{(n)}.\label{mtre}\end{equation}
\par\smallskip
{\bf Step 2.}\quad We prove (\ref{normelp}) and (\ref{aggiunta}) for
any $p\in[1,\infty]$, with the additional assumption that $m=0$ in
$[2,\infty)$. Fix $t>0$. We
consider the function $f$ on
$\R_{+}$ defined by
\[f(\lambda)=e^{\lambda}m(\lambda),\qquad\lambda>0.\]
Then $\|f\|_{(n)}\leq C\|m\|_{(n)}$. Moreover
$M_t=F_t*p_t$, where $F_t$ the kernel of the operator $f(t\LL)$.
So $M_t\in C^{\infty}(G)$ and $X^I M_t=F_t*X^I p_t$. Moreover,
for every $x\in G$ we have
\beq\lefteqn{(1+|x|)^{\alpha}|X^I
M_t(x)|}\\&\leq&C\left(\int_G(1+|y|)^{\alpha}|F_t(y)||X^I
p_t(y^{-1}x)|\,dy+
\int_G|F_t(y)|(1+|y^{-1}x|)^{\alpha}|X^I p_t(y^{-1}x)|\,dy\right).\eeq
Then
\begin{eqnarray}\lefteqn{\|(1+|\cdot|)^{\alpha}X^I M_t(\cdot)\|_{\lp
p}}\nonumber\\&\leq&C\left(\|(1+|\cdot|)^{\alpha}F_t(\cdot)\|_{\lp 1}\|X^I
p_t\|_{\lp{p}}+
\|F_t\|_{\lp{1}}\|(1+|\cdot|)^{\alpha}X^I
p_t(\cdot)\|_{\lp{p}}\right).\label{zzz}\end{eqnarray}
If $|I|=0$ we apply (\ref{ptp}) and
(\ref{mtre}) to both terms of (\ref{zzz}) and we obtain
(\ref{normelp}).
If $0<t\leq 1$ we apply (\ref{xiptp}) and (\ref{mtre}) and we obtain
(\ref{aggiunta}).\par\smallskip
{\bf Step 3.}\quad    We prove (\ref{normelp}) dropping the
additional assumption on $m$. Fix $t>0$ and a non-increasing function
$\varphi\in
C^{\infty}(\R_{+})$ such
that
$\varphi=1$ in $(0,\frac{1}{2})$ and $\varphi=0$ in $[1,\infty)$. Set
\[\psi(\lambda)=\varphi(\frac{\lambda}{2})-\varphi(\lambda),\qquad\lambda>0.\]
So $0\leq\psi\leq 1$ and
supp$\,\psi\subset[\mez,2]$.
Moreover we observe that
\begin{equation}\varphi(\lambda)+\sum_{j=0}^{\infty}\psi(2^{-j}\lambda)=1,\qquad
\lambda>0.
\label{partunita}\end{equation}
So
\[m(\lambda)=\tilde
m(\lambda)+\sum_{j=0}^{\infty}m_j(\lambda),\qquad\lambda>0\]
where
\[\tilde m(\lambda)=m(\lambda)\varphi(\lambda),\qquad\lambda>0\]
and
\[m_j(\lambda)=m(\lambda)\psi(2^{-j}\lambda),\qquad\lambda>0.\]
By (\ref{partunita}) we have
$\sum_{j=0}^{\infty}|m_j|\leq|m|$. So, if
we denote by $\tilde{M}_t$ the kernel of $\tilde m(t\LL)$ and by $M_{j,t}$
the kernel of
$m_j(t\LL)$, by Proposition~\ref{dixmal} we have
\begin{equation}M_t=\tilde{M}_t+\sum_{j=0}^{\infty}M_{j,t}\qquad{\rm in}\
\dig.\label{converg}\end{equation}
We observe that
\[m_j(t\lambda)=h_j(2^{-j}t\lambda),\qquad j\in\N,\ \lambda>0\]
where
\begin{equation}h_j(\sigma)=m(2^j\sigma)\psi(\sigma),\qquad j\in\N,\
\sigma>0.\label{hjl}\end{equation} Since $h_j=0$ in $[2,\infty)$, by the
previous steps
\begin{equation}\|(1+|\cdot|)^{\alpha}M_{j,t}(\cdot)\|_{\lp p}\leq
C(1+2^{-\frac{j}{2}}\sqrt
t)^{\alpha}V(2^{-\frac{j}{2}}\sqrt t)^{-\frac{1}{p'}}\|h_j\|_{(n)}.
\label{mjt}\end{equation}
We observe that the doubling property (\ref{doubling}) implies
\begin{equation}V(2^{-\frac{j}{2}}\sqrt t)^{-\frac{1}{p'}}\leq
K^{\frac{j}{2p'}}V(\sqrt t)^{-\frac{1}{p'}}.\label{vjt}\end{equation}
Moreover we
estimate $\|h_j\|_{(n)}$ by means of (\ref{hjl}): for $r\in\{0,\ldots,n\}$
we have
\begin{equation}h_j^{(r)}(\sigma)=
\sum_{l=0}^r\cb{r}{l}2^{jl}m^{(l)}(2^j\sigma)
\psi^{(r-l)}(\sigma),\qquad j\in\N,\ \sigma>0.\label{hjrl}\end{equation}
Fix an integer $n'>n+\frac{\log_2 K}{2p'}$. By (\ref{hjrl}), if $m\in
C^{n'}(\R_{+})$ with
$\|m\|_{(n')}<\infty$ then
\begin{equation}\|h_j\|_{(n)}\leq C
2^{j(n-n')}\|m\|_{(n')}.\label{fine}\end{equation}
It follows from (\ref{mjt}), (\ref{hjrl}) and (\ref{fine}) that
\begin{equation}\|(1+|\cdot|)^{\alpha}M_{j,t}(\cdot)\|_{\lp p}\leq
C(1+\sqrt t)^{\alpha}V(\sqrt t)^{-\frac{1}{p'}}2^{j(n-n'+\frac{\log_2
K}{2p'})}\|m\|_{(n')}.\label{quasi}\end{equation}
On the other hand
\begin{equation}\|(1+|\cdot|)^{\alpha}\tilde{M}_t(\cdot)\|_{\lp p}\leq
C(1+\sqrt t)^{\alpha}V(\sqrt t)^{-\frac{1}{p'}}\|m\|_{(n)}
\label{tildemt}\end{equation}
since $\|\tilde m\|_{(n)}\leq C\|m\|_{(n)}$.
By (\ref{converg}), (\ref{quasi}) and (\ref{tildemt})
and by the assumption made on $n'$ we obtain
\[\|(1+|\cdot|)^{\alpha}M_t(\cdot)\|_{\lp p}\leq
C(1+\sqrt t)^{\alpha}V(\sqrt t)^{-\frac{1}{p'}}\|m\|_{(n')}.\]
\par\smallskip
{\bf Step 4.}\quad
The proof of (\ref{aggiunta}) without the additional assumption on $m$ is
analogous to Step 3: we
observe that equality (\ref{converg}) implies
\[X^I M_t=X^I\tilde{M}_t+\sum_{j=0}^{\infty}X^I M_{j,t}\qquad{\rm in}\
\dig\]
and then we follow the proof of Step 3: we obtain
\[\|(1+|\cdot|)^{\alpha}X^I M_t(\cdot)\|_{\lp p}\leq
Ct^{-(\frac{d}{2p'}+\frac{|I|}{2})}\|m\|_{(n'')},\qquad 0<t\leq 1\]
where $n''>n+\frac{d}{2p'}+\frac{|I|}{2}$.
\end{proof}\par\smallskip
An immediate
consequence of Proposition~\ref{main} is the following result, which also
generalizes  the
analogous result for stratified groups (see~\cite{Hul}, \cite{Mau}):
\begin{coro} Let $m$ be the restriction on $\R_{+}$ of a function in ${\cal
S}(\R)$. Then the
kernel $M$ of the operator $m(\LL)$ is in
$\sg$.\label{schwartz}\end{coro}

\section{Littlewood--Paley decomposition}\label{Litpal}\leavevmode\par
Fix a non-increasing function $\varphi\in C^{\infty}(\R_{+})$ such that
$\varphi=1$ in $(0,\frac{1}{4})$ and $\varphi=0$ in $[1,\infty)$. Set
\[\psi(\lambda)=\varphi(\frac{\lambda}{4})-\varphi(\lambda),\qquad\lambda>0.\]
So $0\leq\psi\leq 1$ and
supp$\,\psi\subset[\frac{1}{4},4]$.
Moreover we observe that
\begin{equation}\varphi(\lambda)+\sum_{j=0}^N\psi(2^{-2j}\lambda)=\varphi(2^{-2(
N+1)}\lambda),\qquad
N\in\N,\ \lambda>0\label{somma}\end{equation}
and so
\[\varphi(\lambda)+\sum_{j=0}^{\infty}\psi(2^{-2j}\lambda)=1,\qquad\lambda>0.\]
By Corollary \ref{schwartz}, for any $j\in\N$ the kernels of the operators
$S_j=\varphi(2^{-2j}\LL)$ and $\Delta_j=\psi(2^{-2j}\LL)$ are
in $\sg$, so the operators $S_j$ and $\Delta_j$ can be viewed as continuous
operators on $\sig$. By
the spectral theorem, any $f\in\lp 2$ can be decomposed as $f=S_0
f+\sum_{j=0}^{\infty}\Delta_j f$
in
$\lp 2$. The following proposition shows that such decomposition
holds also in $\sg$ and in $\sig$.
\begin{prop} For any $f\in\sg$ and
$u\in\sig$ we have:
\begin{eqnarray}f&=&S_0 f+\sum_{j=0}^{\infty}\Delta_j f\qquad{\rm in}\
\sg;\label{lpf}\\
u&=&S_0 u+\sum_{j=0}^{\infty}\Delta_j u\qquad{\rm in}\
\sig.\label{lpu}\end{eqnarray}\label{lpd}\end{prop}
\begin{proof} We only have to prove (\ref{lpf}), since (\ref{lpu}) follows
by duality. By
(\ref{somma}) we have to prove that
$S_j f\rightarrow f$ in
$\sg$ for
$j\rightarrow \infty$. So we fix
$\alpha\in\N$ and $I\in\bigcup_{\beta\in\N}\{1,\ldots,k\}^{\beta}$ and we
want to prove that
$p_{\alpha,I}(f-S_j f)\rightarrow 0$ for $j\rightarrow
\infty$. Put $N=\max\{n,1+\frac{|I|}{2}\}$, where $n$ is the integer which
appears in Proposition
\ref{main} in the case $p=1$. Then
\begin{equation}f-S_j f=2^{-2jN}m(2^{-2j}\LL)\LL^N f,\qquad
j\in\N\label{fsjf}\end{equation}
where
\begin{equation}m(\lambda)=\frac{1-\varphi(\lambda)}{\lambda^N},\qquad\lambda>0.
\label{mln}\end{equation}
Let $M_j$ be the kernel of the operator $m(2^{-2j}\LL)$.
Since by (\ref{mln}) $m\in C^{\infty}(\R_{+})$ and $\|m\|_{(N)}<\infty$,
Proposition
\ref{main} gives
\begin{equation}\|(1+|\cdot|)^{\alpha}X^I M_j(\cdot)\|_{\lp 1}\leq C
2^{j|I|},\qquad
j\in\N\label{ximj}\end{equation}
where $C$ does not depend on $j$. Then, by (\ref{fsjf}) and
(\ref{ximj}) and since $f\in\sg$, for every $x\in G$ we have
\beq\lefteqn{(1+|x|)^{\alpha}|X^I(f-S_j
f)(x)|}\\&=&2^{-2jN}(1+|x|)^{\alpha}|(\LL^N f*X^I M_j)(x)|\\
&\leq&2^{-2jN}\left(\int_G (1+|y|)^{\alpha}|\LL^N f(y)||X^I
M_j(y^{-1}x)|\,dy+\int_G|\LL^N f(y)|(1+|y^{-1}x|)^{\alpha}|X^I
M_j(y^{-1}x)|\,dy\right)\\ &\leq&C' 2^{j(|I|-2N)}\eeq
where $C'$ does not depend on $j$ or $x$. Since $2N>|I|$, we have that
$p_{\alpha,I}(f-S_j f)\rightarrow 0$ for $j\rightarrow \infty$.
\end{proof}\par\smallskip
{\em Remark:}\quad For $1\leq p<\infty$, since $\sg$ is dense in $\lp p$
and the operators
$\Delta_j$ are uniformly bounded on $\lp p$ by Proposition \ref{main}, it
follows from
(\ref{lpf}) that any
$f\in\lp p$ can be decomposed as $f=S_0 f+\sum_{j=0}^{\infty}\Delta_j f$
in
$\lp p$.\par\medskip
One gets from Proposition \ref{main} an extension
of the classical Bernstein's inequalities for the dyadic blocks (see e.g.
\cite{Peetre},
\cite{Tri1}; see also
\cite[Proposition 3.2]{Lem-book}). In what follows, $d$ is the local
dimension of the group
introduced in Section \ref{Notpre}.
\begin{prop}
For
$I\in\bigcup_{\beta\in\N}\{1,\ldots,k\}^{\beta}$, $1\leq p\leq q\leq
\infty$, $j\in\N$ and $u\in\SS'(G)$ we have:
\beq\|X^I(\sqrt\LL)^{\sigma}S_ju\|_{\lp q}&\leq&C 2^{j(|I|+\sigma+d(\frac
1p -\frac
1q))}
\|S_ju\|_{\lp p},\quad\sigma\geq 0,\\
\|X^I(\sqrt\LL)^{\sigma}\Delta_ju\|_{\lp q}&\leq&C 2^{j(|I|+\sigma+d(\frac
1p -\frac
1q))}
\|\Delta_ju\|_{\lp p},\quad\sigma\in\R,\eeq
where $C$ is a positive constant which depends only on $I,p,q,\sigma$ but
not on $j$ or $u$.
\label{bernstein}\end{prop}
\begin{proof}
We can consider the functions $\tilde \varphi(\lambda)= \varphi
(\frac\lambda4)$ and $\tilde \psi(\lambda)= \varphi(\frac \lambda {16})
-\varphi(4\lambda)$, so that $\tilde \varphi (\lambda) \varphi (\lambda) =
\varphi (\lambda)$
and $ \tilde \psi (\lambda) \psi (\lambda) = \psi (\lambda)$  for all
$\lambda \gt 0$.
The proof of the proposition follows therefore from Young's inequality,
Proposition \ref{main} and
the identities:
\[
\begin{aligned}
S_j u &= \tilde \varphi(2^{-2j}\LL) S_ju;\\
\Delta_ju &= \tilde \psi (2^{-2j} \LL) \Delta_ju.
\end{aligned}
\]\end{proof}

We are now in position to set out the \LP theorem related to
the decomposition (\ref{lpu}):

\begin{theo} Let $1<p<\infty$ and $u\in\SS'(G)$. Then $u\in L^p(G) $ if and
only
if $S_0 u\in L^p(G) $ and $\left( \sum_{j= 0}^\infty
|\Delta_ju|^2\right)^{\frac 12} \in
L^p(G) $. Moreover there exists a constant $C_p>1$, which depends only on
$p$, such that
\[
C_p^{-1}\|u\|_{L^p (G)}\leq\|S_0 u\|_{L^p(G)}+\left \|
\left( \sum_{j= 0}^\infty  |\Delta_ju|^2\right )^{\frac 12}
\right \|_{L^p(G)}\leq
C_p\|u\|_{L^p (G)},\qquad u\in L^p(G).\]\label{llpp}\end{theo}
\begin{proof} Once we have Theorem
\ref{alexop} and equality (\ref{lpu}), the proof of the
Littlewood--Paley theorem in $\R^n$ given for instance in
\cite{Ste} works also in our case.\end{proof}

\section{Besov spaces}\label{Besov}\leavevmode\par
Besov spaces $B^{s,q}_p(G)$ for sub-Laplacians in Lie groups were
studied by many authors.
The usual definition (\cite{Fol}, \cite{Saka})
has been given by means of the heat kernel associated to the sub-Laplacian.
Though in \cite{Saka},  in the case of stratified groups, such spaces are
defined with
$s\in\R$, most applications (see e.g. \cite{Cou-SC}, \cite{SC},
\cite{MouVar}) have concerned essentially the case $s\gt 0$  where
$B^{s,q}_p(G)\subset L^p(G)$. In the Euclidean case there are many
equivalent characterizations
of Besov spaces: a useful reference is given by the books by Triebel
\cite{Tri1}, \cite{Tri2}.
In particular,
the characterizations by the atomic decomposition in the space variable
(\cite{FJ}) and by the dyadic decomposition of the frequency space of the
Fourier transform
(\cite{Peetre}) are often used in the applications. In the case of
unimodular Lie groups,
Skrzypczak in
\cite{Skr2} has given an atomic characterization of Besov spaces with
$s\in\R$. In a more abstract setting, Gal\'e in \cite[Hip\'otesis (${\cal
R}_{\alpha}$)]{Gale}
has determined a sufficient condition for a positive self-adjoint operator
on a Hilbert
space which allows to characterize the corresponding Besov spaces with
$s\gt0$ by a dyadic
decomposition on the spectrum of the operator. In our setting, Gal\'e's
condition is equivalent to
the uniform boundedness of the norms of the operators $\Delta_j
:L^p(G)\freccia L^p(G)$, which
amounts to the uniform estimate  of the norm in $L^1(G)$ of the convolution
kernel of
$\Delta_j$. This is precisely what was pointed out  in the remark after
Proposition
\ref{main}.\par
In this section, we
define on a Lie group of polynomial growth Besov spaces with
$s\in\R$
associated to any sub-Laplacian by means of the \LP
decomposition obtained in Section \ref{Litpal}.
\begin{defi}
Let $s \in \mathbb R$, $1\leq p,q \leq \infty$. We define
\begin{equation*}
    B^{s,q}_{p}(G)=\big\{u \in \sig:\|u\|_{B^{s,q}_{p}(G)}=
    \|S_0u\|_{\lp p} + \left(\sum_{j= 0}^\infty (2^{js}\|\Delta_ju\|_{\lp
p})^q\right)^{1/q}
    < \infty\big\}
\end{equation*}
with obvious modifications in the case $q=\infty$.
\end{defi}

This definition is equivalent to the classical
one through the heat kernel. In fact, once we have Proposition \ref{main}
and the decomposition
(\ref{lpu}), we can repeat the proof given for instance in
\cite[Theorem 5.3]{Lem-book} in the Euclidean setting and we obtain the
following
\begin{prop}\label{besov-eq}
For $s \in \mathbb R,\,1\leq p,q \leq \infty, \, m \geq 0$ such that $m >
s$ and $u\in \sig $, the
following assertions are equivalent:
\begin{itemize}
  \item[i)]  $u \in B^{s,q}_{p}(G);$
  \item[ii)] for all $t>0$, $e^{-t\LL}u\in L^p(G)$  and
             $\int_0^1(t^{-s/2}\|(t\LL)^{m/2} e^{-t\LL}u\|_{\lp p})^q
              \frac{dt}{t} < \infty$.
\end{itemize}
Moreover, the norms\, $\|e^{-\LL}u\|_{\lp p}+(\int_0^1(t^{-s/2}\|(t\LL)^{m/2}
e^{-t\LL}u\|_{\lp p})^q\frac{dt}{t})^{1/q}$ and $\|u\|_{B^{s,q}_p(G)}$ are
equivalent.
\end{prop}\par\smallskip
{\em Remark:}\quad If $s\gt0$, one can replace the condition $e^{-t\LL}u\in
L^p(G)$ by the
equivalent
$u\in L^p(G)$; in fact $B^{s,q}_p(G)\subset B^{0,1}_p(G) \subset \lp p$.
Furthermore, if $s\gt 0$ one can replace the condition
$\int_0^1(t^{-s/2}\|(t\LL)^{m/2}
e^{-t\LL}u\|_{\lp p})^q\frac{dt}{t})\lt \infty$ by the equivalent
$\int_0^\infty(t^{-s/2}\|(t\LL)^{m/2} e^{-t\LL}u\|_{\lp
p})^q\frac{dt}{t})\lt \infty$, due to the
convergence of the integral at infinity.

\small{
\textsc{\ \\
Dipartimento di Ingegneria Gestionale e dell'Informazione, Universit\`a di
Bergamo,\\
Viale Marconi 5, I--24044 Dalmine (BG), Italy} \\
\textit{E-mail:} \texttt{gfurioli@unibg.it}

\vspace{0.1truecm}

\textsc{\ \\
Dipartimento di Scienze Chimiche, Fisiche e Matematiche, Universit\`a
dell'Insubria,\\
Via Valleggio 11, I--22100 Como, Italy} \\
\textit{E-mail:} \texttt{melzi@uninsubria.it}

\vspace{0.1truecm}

\textsc{\ \\
Dipartimento di Matematica, Universit\`a di Genova,\\
Via Dodecaneso 35, I--16146 Genova, Italy} \\
\textit{E-mail:} \texttt{veneruso@dima.unige.it}
}


\begin{thebibliography}{MRS1}

\bibitem[A1]{Ale1}
\textsc{Alexopoulos, G.}, Oscillating multipliers on Lie groups and
Riemannian
manifolds, \emph{T\^ohoku Math. J.} \textbf{46} (1994), 457--468.

\bibitem[A2]{Ale2}
\textsc{Alexopoulos, G.}, Spectral multipliers on Lie groups of polynomial
growth,
\emph{Proc. Amer. Math. Soc.} \textbf{120} (1994), 973--979.

\bibitem[C]{Chr}
\textsc{Christ, M.}, $L^p$ bounds for spectral multipliers on nilpotent
groups,
\emph{Trans. Amer. Math. Soc.} \textbf{328} (1991), 73--81.

\bibitem[CS]{Cou-SC}
\textsc{Coulhon, T.} and \textsc{Saloff-Coste, L.}, Semi-groupes
d'op\'erateurs et
espaces
fonctionnels sur les groupes de Lie, \emph{J. Approx. Theory}
\textbf{65}
  (1991), 176--199.

\bibitem[DM]{DeMa}
\textsc{De Michele, L.} and \textsc{Mauceri, G.}, $H^p$ multipliers on
stratified
  groups, \emph{Ann. Mat. Pura Appl.} \textbf{148} (1987), 353--366.

\bibitem[DiM]{DiMa}
\textsc{Dixmier, J.} and \textsc{Malliavin, P.}, Factorisations de
fonctions et de
vecteurs
  ind\'efiniment diff\'erentiables, \emph{Bull. Sci. Math.} \textbf{102}
(1978),
  305--330.

\bibitem[E]{Ehr}
\textsc{Ehrenpreis, L.}, Some properties of distributions on Lie groups,
\emph{Pacific
  J. Math.} \textbf{6} (1956), 591--605.

\bibitem[F]{Fol}
\textsc{Folland, G. B.}, Subelliptic estimates and function spaces on nilpotent
  Lie groups, \emph{Ark. Mat.} \textbf{13} (1975), 161--207.

\bibitem[FS]{FoSt}
\textsc{Folland, G. B.} and \textsc{Stein, E. M.}, \emph{Hardy spaces on
homogeneous groups},
Math. Notes
\textbf{28}, Princeton Univ. Press, Princeton, 1982.

\bibitem[FJ]{FJ}
\textsc{Frazier, M.} and \textsc{Jawerth, B.}, Decomposition of Besov
spaces, \emph{Indiana
  Univ. Math. J.} \textbf{34} (1985), 777--799.

\bibitem[G]{Gale}
\textsc{Gal{\'e}, J. E.}, Sobre espacios de Besov definidos por medias de
Riesz, in
\emph{Margarita Mathematica en memoria de Jos\'e Javier (Chicho) Guadalupe
  Hern\'andez}, pp.~235--246, Univ. La Rioja, Logro\~{n}o, 2001.

\bibitem[H]{Hor}
\textsc{H{\"o}rmander, L.}, Estimates for translation invariant operators in
  $L^p$ spaces, \emph{Acta Math.} \textbf{104} (1960), 93--140.

\bibitem[Hu]{Hul}
\textsc{Hulanicki, A.}, A functional calculus for Rockland operators on
nilpotent
  Lie groups, \emph{Studia Math.} \textbf{78} (1984), 253--266.

\bibitem[L]{Lem-book}
\textsc{Lemari{\'e}-Rieusset, P. G.}, \emph{Recent developments in the
  Navier--Stokes problem}, Chapman \& Hall /CRC Research Notes in
  Math. \textbf{431}, Chapman \& Hall /CRC, Boca Raton, 2002.

\bibitem[M]{Mau}
\textsc{Mauceri, G.}, Maximal operators and Riesz means on stratified
groups, in
\emph{Sympos. Math.} \textbf{29}, pp.~47--62, Academic Press, New York, 1987.

\bibitem[MM]{MaMe}
\textsc{Mauceri, G.} and \textsc{Meda, S.}, Vector-valued multipliers on
stratified
groups,
\emph{Rev. Mat. Iberoamericana} \textbf{6} (1990), 141--154.

\bibitem[Me]{Mel}
\textsc{Melrose, R.}, Propagation for the wave group of a positive subelliptic
  second-order differential operator, in \emph{Hyperbolic equations and related
topics
  (Katata/Kyoto 1984)}, pp.~181--192, Academic Press, Boston, 1986.

\bibitem[MRS1]{MRS1}
\textsc{M{\"u}ller, D.}, \textsc{Ricci, F.} and \textsc{Stein, E. M.},
Marcinkiewicz
multipliers and
  multi-parameter structure on Heisenberg (-type) groups, I, \emph{Invent.
Math.}
  \textbf{119} (1995), 199--233.

\bibitem[MRS2]{MRS2}
\textsc{M{\"u}ller, D.}, \textsc{Ricci, F.} and \textsc{Stein, E. M.},
Marcinkiewicz
multipliers and
  multi-parameter structure on Heisenberg (-type) groups, II, \emph{Math.
Z.} \textbf{221} (1996), 267--291.

\bibitem[MV]{MouVar}
\textsc{Mustapha, S.} and \textsc{Varopoulos, N. Th.}, Comparaison
H{\"o}lderienne des
  distances sous-elliptiques et calcul $S(m,g)$, \emph{Potential Anal.}
\textbf{4}
  (1995), 415--428.

\bibitem[P]{Peetre}
\textsc{Peetre, J.}, \emph{New thoughts on Besov spaces}, Duke Univ. Math.
Series \textbf{1},
Duke
  University, Durham, 1976.

\bibitem[S]{Saka}
\textsc{Saka, K.}, Besov spaces and Sobolev spaces on a nilpotent Lie group,
\emph{T\^ohoku Math. J.} \textbf{31} (1979), 383--437.

\bibitem[Sa]{SC}
\textsc{Saloff-Coste, L.}, Analyse sur les groupes de Lie \`a croissance
  polyn\^omiale, \emph{Ark. Mat.} \textbf{28} (1990), 315--331.

\bibitem[S1]{Skr1}
\textsc{Skrzypczak, L.}, Atomic decompositions on manifolds with bounded
geometry,
\emph{Forum Math.} \textbf{10} (1998), 19--38.

\bibitem[S2]{Skr2}
\textsc{Skrzypczak, L.}, Besov spaces and Hausdorff dimension for some
  Carnot--Carath\'eodory metric spaces, \emph{Canad. J. Math.} \textbf{54}
  (2002), 1280--1304.

\bibitem[St]{Ste}
\textsc{Stein, E. M.}, \emph{Harmonic analysis}, Princeton Math. Series
\textbf{43}, Princeton Univ.
Press, Princeton, 1993.

\bibitem[T1]{Tri1}
\textsc{Triebel, H.}, \emph{Theory of function spaces}, Monogr. in
Math. \textbf{78}, Birkh\"auser Verlag, Basel, 1983.

\bibitem[T2]{Tri2}
\textsc{Triebel, H.}, \emph{Theory of function spaces II}, Monogr. in
Math. \textbf{84}, Birkh\"auser Verlag, Basel, 1992.

\bibitem[VSC]{VSC}
\textsc{Varopoulos, N. Th.}, \textsc{Saloff-Coste, L.} and \textsc{Coulhon,
Th.}, \emph{Analysis
and
  geometry on groups}, Cambridge Tracts in Math. \textbf{100}, Cambridge
Univ. Press, Cambridge,
1992.

\end{thebibliography}
\end{document}